\documentclass[10pt]{article}
\pdfoutput=1
\usepackage{etoolbox}

\newtoggle{fullversion}
\toggletrue{fullversion}
\newtoggle{ovcompile}
\togglefalse{ovcompile}

\usepackage{times}

\iftoggle{ovcompile}
{
	
}{}

\usepackage[utf8]{inputenc}
\usepackage[bookmarks=true]{hyperref}
\usepackage{xcolor}
\usepackage{etoolbox}
\usepackage{upgreek}
\usepackage{amsmath}
\usepackage{comment,cleveref}

\usepackage{mathtools}
\usepackage{graphicx}
\usepackage{dsfont} 

\usepackage{amssymb}
\usepackage{mathtools}

\iftoggle{fullversion}
{\newcommand{\cdcstyle}{}
\newcommand{\cdcpar}[1]{\paragraph{#1}}
\newcommand{\cdcparsssec}[1]{\subsubsection{#1}}
\usepackage{fullpage}
\newcommand{\fullvtext}[1]{#1}
}
{\newcommand{\cdcstyle}{\textstyle}
\newcommand{\cdcpar}[1]{\noindent\textbf{#1:}}
\newcommand{\cdcparsssec}[1]{\noindent\textbf{#1:}}
\newcommand{\fullvtext}[1]{}
}

\Crefname{equation}{Eq.}{Eqs.}
\Crefname{assumption}{Assumption}{Assumptions}
\Crefname{proposition}{Proposition}{Props}

\newcommand{\set}[1]{\left\{ #1\right\}}

\iftoggle{fullversion}
{
	\usepackage{amsthm}
	\newtheorem{theorem}{Theorem}
	\newtheorem{assumption}{Assumption}
	\newtheorem{definition}{Definition}
	\newtheorem{lemma}{Lemma}
	\newtheorem{remark}{Remark}
	\newtheorem{proposition}{Proposition}

}
{
	\newtheorem{assumption}{Assumption}
	\newtheorem{lemma}{Lemma}
	\newtheorem{theorem}{Theorem}
	
	\newtheorem{definition}{Definition}
	\newtheorem{remark}{Remark}
	\newtheorem{proposition}{Proposition}
}

\newcommand{\R}{\mathbb{R}}

\newcommand{\N}{\mathbb{N}}

\newlength\tindent
\let\sp=\sqparen

\renewcommand{\indent}{\hspace*{\tindent}}
\newcommand{\twnote}[1]%
    {\textcolor{orange}{\textbf{TW: #1}}}
\newcommand{\msnote}[1]%
    {\textcolor{cyan}{\textbf{MS: #1}}}

\newcommand{\rmd}{\mathrm{d}}
\newcommand{\dt}{\mathrm{d}t}
\newcommand{\FOS}{\mathtt{FOS}}
\newcommand{\calL}{\mathcal{L}}

\newcommand{\nablatwo}{\nabla^{\,2}}

\newcommand{\Jlin}{J^{\mathtt{jac}}}
\newcommand{\Ltwo}{\calL^2}

\newcommand{\xbar}{\bar{x}}
\newcommand{\ubar}{\bar{u}}
\newcommand{\Vtil}{\tilde{\mathcal{V}}}
\newcommand{\Vjacst}{\mathcal{V}^{\mathtt{jac},\star}}
\newcommand{\ddt}{\frac{\rmd}{\rmd t}}
\newcommand{\Atil}{\tilde{A}}
\newcommand{\Btil}{\tilde{B}}

\newcommand{\xtil}{\tilde{x}}
\newcommand{\util}{\tilde{u}}
\newcommand{\xcheck}{\check{x}}
\newcommand{\ucheck}{\check{u}}
\newcommand{\calC}{\mathcal{C}}
\newcommand{\hatu}{\bar{u}}

\title{On the Stability of Nonlinear Receding Horizon Control: \\A Geometric Perspective}
\author{Tyler Westenbroek$^{*1}$, Max Simchowitz$^{*1}$, Michael I. Jordan$^{1,2}$, S. Shankar Sastry$^{1}$
\thanks{$^*$ Indicates equal contribution.}
\thanks{$^{1}$Department of Electrical Engineering and Computer Sciences, University of California at Berkeley, USA.}%
\thanks{$^2$Department of Statistics, University of California at Berkeley, USA.}}

\begin{document}

\maketitle
\begin{abstract}

The widespread adoption of nonlinear Receding Horizon Control (RHC) strategies by industry has led to more than 30 years of intense research efforts to provide stability guarantees for these methods. However, current theoretical guarantees require that each (generally nonconvex) planning problem can be solved to (approximate) global optimality, which is an unrealistic requirement for the derivative-based local optimization methods generally used in practical implementations of RHC. This paper takes the first step towards understanding stability guarantees for nonlinear RHC when the inner planning problem is solved to first-order stationary points, but not necessarily global optima. Special attention is given to feedback linearizable systems, and a mixture of positive and negative results are provided.
We establish that, under certain strong conditions, first-order solutions to RHC exponentially stabilize linearizable systems. Surprisingly, these conditions can hold even in situations where there may be \emph{spurious local minima.}  Crucially, this guarantee requires that state costs applied to the planning problems are in a certain sense `compatible' with the global geometry of the system, and a simple counter-example demonstrates the necessity of this condition. These results highlight the need to rethink the role of global geometry in the context of optimization-based control. 

\end{abstract}

\section{Introduction}\label{sec:intro}

The global stabilization of nonlinear systems is one of the most fundamental and challenging problems in control theory. In principle, the search for a stabilizing controller can be reduced to finding a control Lyapunov function \cite{SONTAG1989117}, yet the search for such a function may be just as challenging as the search for a controller. Historically, the framing of these equivalent synthesis problems has been guided by two distinct perspectives: geometric control and optimal control.

Geometric control is a broad term that loosely refers to design methodologies which systematically exploit global system structure to achieve a control objective for a specific class of systems \cite{isidori1990output,kokotovic1992joy, sastry2013nonlinear}. In the context of stabilization, this has led to constructive procedures for synthesizing stabilizing controllers or control Lyapunov functions for many important classes of systems (feedback linearizable, strict feedback, etc). Despite the broad impact of this collection of techniques, the traditional criticism of these approaches is that they often require significant human ingenuity and system-specific analysis to implement.

Conversely, optimal control offers a way to sidestep these challenges by formulating certain infinite-horizon problems which can be solved to synthesize optimal stabilizing controllers automatically.  Unfortunately, exact solutions to the infinite-horizon problem require solving the Hamilton-Jacobi-Bellman partial differential equation, whose exact solution is generally impractical in state spaces of large or even modest dimensionality \cite{bellman1966dynamic}. The response to these limitations has been the rise of \emph{receding horizon} or \emph{model-predictive} control strategies  \cite{mayne2000constrained,jadbabaie2001stability,primbs1999nonlinear}, which attempt to stabilize the system by solving a sequence of more tractable open-loop optimal control problems along the system trajectory. 

Due to the widespread adoption of these methods by practitioners, providing stability guarantees for nonlinear receding horizon methods has been an active area of research for theoreticians for more than 30 years. However, there is a distinct gap between the theoretical models of RHC studied in the literature and the practical implementations favored in industry. Specifically, while current theoretical results require that each optimal control problem posed in the RHC scheme can be solved to (near) global optimality, practical implementations of RHC, which usually employ local derivative-based search algorithms, can only be guaranteed to find approximate stationary points of each planning problem. Thus, since nonlinear RHC planning problems are generally nonconvex, providing stability guarantees for derivative-based nonlinear RHC implementations remains a largely open problem. Due to this gap, designing cost functions for derivative-based RHC schemes which reliably guide the local search algorithms towards stabilizing solutions remains a heuristic-driven and often time-consuming process. 

In this paper we aim to link the geometric and optimal control perspectives, demonstrating how such a link provides new insights into the stability of derivative-based nonlinear RHC implementations. Specifically, by studying how the RHC cost functions interact with both the local and global geometry of the control system, we illustrate how the choice of cost function can lead either to provably stable behavior, or to failure modes where RHC control schemes get `stuck' at undesirable stationary points. Our negative results are related through two counter-examples while our positive results are given in Theorem 1, which provides sufficient conditions which ensure that all (approximate) first-order stationary points of the RHC optimal control problems correspond to open-loop state trajectories which decay exponentially to the origin. We use this result to provide stability guarantees for nonlinear RHC when the implementation relies on derivative-based descent methods (\Cref{thm:RHC}), provided that the RHC planning horizon is of sufficient (though modest) length.

Our sufficient conditions for exponential stability informally require $(1)$ the state and input costs are strongly convex,  $(2)$ the Jacobian linearization of the dynamics along every trajectory is uniformly stabilizable, $(3)$ the control cost is sufficiently small when compared to the state cost, $(4)$ the (time-varying, affine) first-order expansion of the system dynamics along each trajectory satisfies a local `matching condition', meaning that the affine term can be cancelled out by an appropriate choice of input, and $(5)$ the nonlinearity in the input is sufficiently small. We note that assumption $(1)$ is natural in practical implementations of RHC, and some (possibly weaker) stabilizability condition as in $(2)$ is clearly necessary to ensure exponential stability. It is unclear whether $(5)$ is necessary, and we leave its study to future work. We note that conditions $(1)-(2)$ and $(4)$ may be satisfied in some coordinate systems for the state but not others. Thus applying our sufficient conditions (or using them to design a `good' RHC cost functional) may require finding an appropriate coordinate system for the system.

To shed light on conditions $(3)$ and $(4)$, and to establish their necessity, we examine a class of feedback linearizable systems which satisfy conditions $(2)$, $(4)$ and $(5)$ in an appropriate choice of coordinates. For this class of systems conditions $(1)$ and $(4)$ require that the state costs are strongly convex \emph{in the linearizing coordinates}. To demonstrate the necessity of this strong geometric condition, we examine a model for a flexible-joint manipulator which is full-state linearizable. A natural state cost is designed which is convex in the `original' non-linearizing coordinates (where condition $(4)$ is violated), but analysis reveals that the cost is non-convex in the linearizing coordinates (where condition $(4)$ is satisfied). Thus, the chosen cost function is in some sense `incompatible' with the geometry of the system as we are forced to pick a coordinate system in which either $(1)$ or $(4)$ are violated. Due to this mismatch, we are able to identify initial conditions from which derivative-based RHC schemes will fail to stabilize the system and get stuck at undesirable stationary points. To address condition $(3)$ we also investigate a model for the simple inverted pendulum, which is in the class of linearizable systems discussed above. We demonstrate that even when the other four conditions are satisfied RHC may again fail to stabilize the system if the penalty on the input is too large.

 Unlike prior works \cite{jadbabaie2001stability} which require global optimality for each RHC planning problem, the stability issues for derivative-based RHC implementations discussed above cannot be overcome by simply increasing the prediction horizon.
 Indeed, our counterexamples show that conditions $(3)$ and $(4)$ are necessary even when arbitrarily long planning horizons are used. The key difference here is that derivative-based planners can converge to overly-myopic sequences of control inputs, even with long planning horizons, due to the myopic nature of the optimization landscape (i.e., the presence of local minima).

In sum, the results of this paper indicate that the stability of derivative-based nonlinear RHC schemes may be fragile unless the interaction between the geometry of the control system and cost functions are carefully considered. Fortunately, our positive results indicate that concepts from the geometric control literature may provide constructive techniques for designing RHC cost functionals which provably guide local search algorithms towards stabilizing solutions.

 \fullvtext{In addition, we establish that there can be stabilizing trajectories which satisfy our conditions, but are not global optima; therefore, the results we establish hold for cost functionals which need not exhibit quasi-convexity or any other obvious form of ``hidden'' convexity.}

\cdcpar{Further Background on RHC}
For basic background on RHC, we refer the reader to any of a number of comprehensive reviews on RHC (see, e.g., \cite{mayne2000constrained,morari1999model}).  Simplifying considerably, previous theoretical nonlinear RHC formulations fall into either \emph{constrained} approaches and \emph{unconstrained} approaches. Constrained RHC formulations directly enforce stability by either constraining the terminal predictive state to lie at the origin \cite{mayne1988receding}, or using inequality constraints to force the system into a neighborhood containing the origin, and then stabilizing the system using a local controller \cite{michalska1993robust}. The usual critique of these methods \cite{jadbabaie2001unconstrained} is that the satisfaction of the relevant constraints may be overly demanding computationally in an online implementation. In contrast, unconstrained approaches implicitly enforce stability by either using an appropriate CLF as the terminal cost \cite{jadbabaie2001unconstrained} or a sufficiently long prediction horizon \cite{jadbabaie2001stability}. As alluded to above, most of these stability guarantees require that a globally optimal solution can be found for each prediction problem. Several approaches provide stability guarantees using sub-optimal solutions, but generally require that an initial feasible solution is available \cite{scokaert1999suboptimal}, which may be restrictive in high-performance real-time scenarios, or require the availability of a CLF ~\cite{primbs1999nonlinear, jadbabaie2001unconstrained}, which implies that the stabilization problem has already been solved. Thus, in this paper we study unconstrained RHC formulations which use general terminal costs, and aim to provide stability guarantees which only require that a stationary point of each optimization problem can be found. We feel that this accurately reflects the spirit of optimization-based control---to stabilize the system with minimal system-specific knowledge---as well as the practical computational constraints facing practitioners. 
\iftoggle{fullversion}
{
	
}
{

\noindent\textbf{Remark On Proofs:}  The proofs of \Cref{lemma:approx_prob,lem:underbound}, \Cref{prop:fbl1}, and \Cref{thm:RHC} can be found in \cite{techreport}. 
}

\section{Preliminaries}
\newcommand{\range}{\mathrm{range}}
\newcommand{\tinyddt}{\tfrac{\rmd}{\rmd t}}
\newcommand{\UT}{\mathcal{U}_T}
\newcommand{\qedd}{\scalebox{0.8}{$\blacksquare$}}
\newcommand{\ieeethm}{\itshape}
This paper studies control systems of the form
\begin{equation}\label{eq:dynamics}
    \dot{x}(t)  =   F(x(t),u(t)),
\end{equation}
where $x \in \R^n$ is the state and $u \in \R^m$ the input, and $\dot{x}(t)= \ddt x(t)$ denotes time derivatives. We make the following assumptions about the vector field $F\colon \R^n \times \R^m \to \R^n$:
\begin{assumption}\label{ass:eq_point}
\ieeethm The origin is an equilibrium point of \eqref{eq:dynamics}, namely, $F(0,0)=0$. 
\end{assumption}
\begin{assumption}\label{ass:v_field_regularity} \ieeethm
The vector field $F$ is continuously differentiable. Furthermore, there exist constants $L_F>0$ such that for each $x_1,x_2 \in \R^n$ and $u_1, u_2 \in \R^m$ we have:
\begin{equation*}
\|F(x_1,u_1) - F(x_2,u_2)\|_2 \leq L_F \left( \|x_1 - x_2\|_2 + \| u_1 -u_2\|_2\right).
\end{equation*}
\end{assumption}

Taken together, these standard assumptions support the global existence and uniqueness of solutions to \eqref{eq:dynamics} on compact intervals of time \cite[Proposition 5.6.5]{polak2012optimization}.

 The primary object of study in is paper will be finite horizon cost functionals $J_T(\cdot;x_0) \colon \UT \to \R$ of the form
\begin{align}\nonumber
 J_{T}(\tilde{u};x_0)& = \int_{0}^{T}  Q(\tilde{x}(\tau)) + R(\tilde{u}(\tau)) \dt +V(\tilde{x}(T))\\
& \ \   \text{s.t. } \dot{ \tilde{x}}(t) = F(\tilde{x}(t),\tilde{u}(t)) , ~~ \tilde{x}(0) = x_0, \label{eq:JT_dynamics}
\end{align}
where $T>0$ is a finite prediction horizon, $x_0\in \R^n$ is the initial condition for \eqref{eq:dynamics},  and the space of admissible inputs is given by $\UT \colon = \Ltwo\left( \sp{0,T}, \R^m \right) \cap \calL^{\infty}\left(\sp{0,T}, \R^m\right)$. Here, $Q \colon \R^n \to \R$ is the running cost applied to the state, $R\colon \R^m \to \R$ is the running cost applied to the input, and $V \colon \R^n \to \R$ a penalty for the terminal state. For now we assume that each of these maps is continuously differentiable. 

\newcommand{\delutil}{\updelta \tilde{u}}
\newcommand{\delu}{\updelta u}
\newcommand{\delx}{\updelta x}

\newcommand{\Der}{\mathrm{D}}
\subsection{Linearizations and (Approximate) Stationary Points}\label{sec:math}
Next, we briefly review a few basic facts from the calculus of variations which are essential for understanding out results. We endow $\Ltwo\left( \sp{0,T}, \R^m \right)$ with the usual inner product and  norm, denoted $\langle \cdot, \cdot \rangle \colon \Ltwo(\sp{0,T}, \R^m) \times \Ltwo(\sp{0,T}, \R^m) \to \R$ and $\|\cdot\|_2: \Ltwo(\sp{0,T}, \R^m)  \to \R $. Under Assumptions \ref{ass:v_field_regularity} and \ref{asm:cost_smooth}, directional (Fr\'{e}chet) derivatives of $J_T(\cdot,x_0)$ are guaranteed to exist \cite[Theorem 5.6.8]{polak2012optimization} as there is a well-defined gradient at each point in the optimization space. We denote the directional Fr\'echet derivative of $J_T(\cdot,x_0)$ at the point $\tilde{u} \in \UT$ in the direction $\delu\in \UT$ by $\Der J_T(\tilde{u};x_0; \delu)$. The gradient $\nabla J_T(\tilde{u} ;x_0) \in \Ltwo(\sp{0,T}, \R^m)$ is the unique object satisfying, for each $\delu \in \UT$, 
\begin{align}\label{eq:grad1}
\Der J_T(\tilde{u};x_0; \delu) 
&= \cdcstyle \int_0^T \langle \nabla J_{T} ( \tilde{u}; x_0)(t) , \delu(t) \rangle \dt, \quad 
\end{align}
or more compactly, $\Der J_T(\tilde{u};x_0; \delu) = \langle \nabla J_T(\tilde{u} ;x_0), \delu \rangle $. The following notions from the optimization literature are crucial for understanding our technical results: 

\begin{definition}\ieeethm
We say that an input $\tilde{u}$ is a \textbf{\emph{first-order stationary point}} ($\FOS$) if $\nabla J_{T}(\tilde{u};x_0) = 0$. We say that $\tilde{u}$ is an $\epsilon$-$\FOS$ if $\|\nabla J_{T}(\tilde{u};x_0)\|_2 \leq \epsilon$. 
\end{definition}

In practice, derivative-based descent algorithms take an infinite number of iterations to converge to exact stationary points, thus our analysis will primarily focus on the approximate stationary points of $J_T(\cdot,x_0)$, as these can be reached in a finite number of iterations.

Finally, we discuss how to calculate the gradient $\nabla J_T(\tilde{u} ;x_0)$ using first-order expansions of the system dynamics and cost functions. Let $(\tilde{x}(\cdot),\tilde{u}(\cdot))$ denote the state-input pair of \eqref{eq:dynamics} defined on the interval $\sp{0,T}$ such that $\tilde{x}(0) = x_0$ and define for each $t \in \sp{0,T}$ \begin{equation*}
\tilde{A}(t) = \frac{\partial}{\partial x}F(\tilde{x}(t),\tilde{u}(t)), \ \ \ \tilde{B}(t) = \frac{\partial }{\partial u}F(\tilde{x}(t),\tilde{u}(t)).
\end{equation*}

\begin{definition}\vspace{0.1cm}\ieeethm
Let $(\tilde{x}(\cdot)$,$\tilde{u}(\cdot))$, $\tilde{A}(\cdot)$ and $\tilde{B}(\cdot)$ be defined as above. We refer to the time-varying linear system $(\tilde{A}(\cdot), \tilde{B}(\cdot))$ as the \textbf{\emph{Jacobian linearization}} of the vector field $F$ along the trajectory $(\tilde{x}(\cdot), \tilde{u}(\cdot))$.
\end{definition}

The Jacobian linearization $(\tilde{A}(\cdot),\tilde{B}(\cdot))$ can be used to construct a first-order approximations to trajectories near $(\tilde{x}(\cdot),\tilde{u}(\cdot))$ as follows. Let $\delu \in \UT$ be a small admissible perturbation to the input and let $(\hat{x}(\cdot), \hat{u}(\cdot))$ satisfy \eqref{eq:dynamics} and $\hat{u} = \tilde{u} + \delu$ and $\hat{x}(0) = \tilde{x}(0) = x_0$. Then a first-order approximation to $\hat{x}(\cdot)$ is given by 
\begin{equation}\label{eq:approx_flow}
    \hat{x}(\cdot) \approx \bar{x}(\cdot) := \tilde{x}(\cdot) + \delx(\cdot)
\end{equation}
where $\delx \colon \sp{0,T} \to \R^n$ solves
    $\dot{\delx} (t) = \tilde{A}(t) \delx(t) + \tilde{B}(t)\delu(t)$ with $\delx(0) = 0$,
and the approximation in \eqref{eq:approx_flow} suppresses higher-order terms involving $\delu$. Regarding gradients as row vectors, by \cite[Theorem 5.6.3]{polak2012optimization} we have
\begin{equation}\label{eq:grad_calc}
    \nabla J_T(\tilde{u} ;x_0)(t) = p(t)\tilde{B}(t) + \nabla R(\tilde{u}(t))
\end{equation}
where the \emph{co-state} $p \colon \sp{0,T} \to \R^{1\times n}$ satisfies
\begin{equation}\label{eq:costate}
    -\dot{p}(t) = p(t)\tilde{A}(t) + \nabla Q(\tilde{x}(t)) \ \ \ p(t) = \nabla V(\tilde{x}(T)).
\end{equation}
Thus, equations \eqref{eq:grad_calc} and \eqref{eq:costate} reveal that the gradient of the objective can be efficiently computed using a `backwards pass' along the nominal trajectory $(\tilde{x}(\cdot),\tilde{u}(\cdot))$ and the linearizations of the vector field and costs along this curve.

\newcommand{\deldotx}{\updelta \dot{x}}

\subsection{Convex Time-Varying Approximations to $J_T(\cdot,x_0)$}\label{sec:jac}
Our primary goal throughout the paper is study the properties of (approximate) stationary points of $J_T(\cdot;x_0)$, and our primary analytical tool will be a family of convex approximations constructed using the Jacobian linearization around particular trajectories of the system. To begin constructing these approximations, observe that the evolution of the estimate $\bar{x}(\cdot)$ in \eqref{eq:approx_flow} is given by
\begin{align} 
\dot{\bar{x}}(t) &= \dot{\tilde{x}}(t) + \deldotx(t) \nonumber\\
& = F(\tilde{x}(t), \tilde{u}(t)) + \tilde{A}(t)\delx(t) + \tilde{B}(t)\delu(t)\nonumber\\ 
& = \tilde{A}(t)\bar{x}(t) + \tilde{B}(t)  \bar{u}(t) + \tilde{d}(t), \text{ where } \label{eq:approx_dyn}\\
\tilde{d}(t) &:= F(\tilde{x}(t),\tilde{u}(t)) - \tilde{A}(t)\tilde{x}(t) - \tilde{B}(t)\tilde{u}(t).\label{eq:drift_def}
\end{align}
We use these time-varying dynamics to approximate $J_T(\cdot; x_0)$ near the point $\tilde{u}$ with the cost functional $\Jlin_{T}(\cdot;x_0,\tilde{u}) \colon \UT \to \R$ defined as follows:
\begin{align}\label{eq:jlin}
 & \Jlin_{T}(\bar{u};x_0,\tilde{u}) = {\textstyle\int_{0}^{T}}  Q(\bar{x}(\tau)) + R(\bar{u}(\tau)) \dt +V(\bar{x}(T)) \nonumber\\
&~~ \text{\scalebox{.9}{{s.t. } $\dot{\bar{x}}(t) = \tilde{A}(t)\bar{x}(t) + \tilde{B}(t)\bar{u}(t) + \tilde{d}(t), ~~ \bar{x}(0) = x_0$}}.
\end{align}
 The following result, which follows from a direct comparison of the formulas for the gradients of $J_T(\cdot; x_0)$ and $\Jlin_T(\cdot;x_0, \tilde{u})$ at the point $\tilde{u}$, motives this construction:
\begin{lemma}\label{lemma:approx_prob}\ieeethm
 For any input $\tilde{u}(\cdot) \in \UT$ we have  \iftoggle{fullversion}
 {
 \begin{align*}
 J_T(\tilde{u};x_0) =  \Jlin_T(\tilde{u};x_0, \tilde{u}), \quad \text{and} \quad \nabla J_T(\tilde{u};x_0) =  \nabla\Jlin_T(\tilde{u};x_0, \tilde{u}).
 \end{align*}
 }
 {
  $ J_T(\tilde{u};x_0) =  \Jlin_T(\tilde{u};x_0, \tilde{u})$ and $\nabla J_T(\tilde{u};x_0) =  \nabla\Jlin_T(\tilde{u};x_0, \tilde{u})$. 
 }
\end{lemma}
In particular, this comparison result will allow us to characterize the (approximate) stationary points of the full nonlinear problem using the local approximation, as the stationary points for these two problem agree along the nominal control input. As we shall see, studying the linearized problem will enable us to connect the local (first-order) structure of the problem to the global structure of the control system. 

\section{Sufficient Conditions for Exponentially Decaying First-Order Stationary Points}\label{sec:main_result}
We begin our analysis by providing sufficient conditions which ensure that \emph{all} (approximate) stationary points of $J_T(\cdot;x_0)$ decay exponentially to the origin at a rate that is independent of $T\geq 0$ and $x_0 \in \R^n$. This is a strong condition which will enable us to provide global exponential stability guarantees for RHC schemes which use derivative-based iterative optimization schemes in Section \ref{sec:FORHC}. After stating the main result and its proof, we draw the connection to feedback linearization and investigate examples which highlight the necessity of some of our stronger assumptions.

\subsection{Sufficient Conditions for Exponentially Decaying (Approximate) First-Order Stationary Points}\label{subsec:exp_convergence}
\newcommand{\ulin}{u^{\mathtt{lin}}}

\fullvtext{In addition, we provide a further counterexample demonstrating that our assumptions permit stationary points which are not globally optimal.} 
We begin by introducing the Assumptions required for the proof of Theorem \ref{thm:plan}. We emphasize some of these assumptions are highly dependent on the particular set of coordinates chosen for the state, and may be satisfied in certain coordinate systems but not others. Thus, applying our sufficient conditions requires finding a coordinate system in which the following conditions hold. Later we relate these conditions to the coordinate systems that arise when performing feedback linearization, which succinctly capture the underlying geometry of the control system.

Our first assumption is strong convexity and smoothness of the running and terminal cost functions: 

\begin{assumption}\label{asm:cost_smooth} \ieeethm We assume that $Q(\cdot),R(\cdot),V(\cdot)$ are twice-continuously differentiable 
functions, with $Q(0) =R(0)=V(0) = 0$,
whose Hessians satisfy the pointwise bounds $\alpha_Q I \preceq \nablatwo Q \preceq \beta_Q I$, $\alpha_R I \preceq \nablatwo R \preceq \beta_R I $, and $ \alpha_V I \preceq \nablatwo V \preceq \beta_V I $ for constants $0 < \alpha_Q  \le \beta_Q$, $0 < \alpha_R \le \beta_R$, and $0\leq \alpha_V \leq \beta_V $. 
\end{assumption}

Note that when Assumption \ref{asm:cost_smooth} is satisfied the optimization in \eqref{eq:jlin} is strongly convex. Thus, by Lemma \ref{lemma:approx_prob}, $\tilde{u}$ is a stationary point of $J_T(\cdot,x_0)$ if and only if it is the \emph{global minimizer} of $\Jlin_{T}(\cdot;x_0,\tilde{u})$. Due to the convexity of the approximation, it is much easier to study the properties of stationary points of $J_T(\cdot,x_0)$ using $\Jlin_{T}(\cdot;x_0,\tilde{u})$ rather than the original functional. This observation is a key insight in our poof technique. The following result extends the above discussion to approximate stationary points of $J_T(\cdot,x_0)$: 
\begin{lemma}\label{lem:comparison} (Approximate $\FOS$) \ieeethm Suppose \Cref{asm:cost_smooth} holds, and that $\tilde{u}$ is an $\epsilon$-$\FOS$ of $J_T(\cdot;x_0)$. Then, 
\begin{align*}
 J_T(\tilde{u};x_0) \le \min_{u}\Jlin_{T}(u;x_0,\tilde{u}) + \tfrac{\epsilon^2}{2\alpha_R}.
\end{align*}
\end{lemma}
\begin{proof} \Cref{asm:cost_smooth} implies that  $\Jlin_T(\bar{u};x_0,\tilde{u}) - \alpha_R \|\bar{u}\|^2$ is convex, and thus the Polyak-\L{}ojasiewicz inequality  holds: $\Jlin_T(\bar{u};x_0,\tilde{u}) \le  \min_{u}\Jlin_{T}(u;x_0,\tilde{u}) +  \|\nabla \Jlin_T(\bar{u};x_0,\tilde{u})\|_2^2/2\alpha_R$ .   Since $\tilde{u}$ is an $\epsilon$-$\FOS$ of $J_T(\cdot;x_0)$ and 
 $ \nabla J_T(\tilde{u};x_0) = \nabla \Jlin_{T}(\tilde{u};x_0,\tilde{u})$, it follows that $\Jlin_T(\tilde{u};x_0,\tilde{u}) \le  \min_{u}\Jlin_{T}(u;x_0,\tilde{u}) +  \epsilon^2/2\alpha_R$. Using $ J_T(\tilde{u};x_0) = \Jlin_{T}(\tilde{u};x_0,\tilde{u})$ concludes.
\end{proof}

Next we place restrictions on the local structure of the control system along each of its trajectories:

\begin{assumption}\label{asm:noise_structure}\ieeethm
Along each system trajectory $(\tilde{x}(\cdot), \tilde{u}(\cdot))$ the drift term in \eqref{eq:drift_def} satisfies $\tilde{d}(t) \in \range(\tilde{B}(t))$. 
\end{assumption}

\begin{assumption}\label{asm:stabilizable} \ieeethm
There exists $\gamma>0$ such that for each time horizon $T\geq0 $, $x_0 \in \R^n$ and system trajectory $(\tilde{x}(\cdot),\tilde{u}(\cdot))$ of length $T$ with $\tilde{x}(0) = x_0$ we have
\begin{equation}
    \inf_{\hat{u}(\cdot)} \cdcstyle \int_{0}^{T}\| \hat{x}(t)\|_2^2 + \|\hat{u}(t)\|_2^2\dt + \|\hat{x}(T)\|_2^2 \leq \gamma\| x_0\|_2^2,
\end{equation}
where $\dot{\hat{x}}(t) = \tilde{A}(t)\hat{x}(t) + \tilde{B}(t)\hat{u}(t)$ and $\hat{x}(s) =x_0$ and $(\tilde{A}(\cdot),\tilde{B}(\cdot))$ is the Jacobian linearization along $(\tilde{x}(\cdot),\tilde{u}(\cdot))$. 
\end{assumption}

 In the language of nonlinear control theory, Assumption \ref{asm:noise_structure} is known as a \emph{matching condition} \cite[Chapter 9.4]{sastry2013nonlinear}. The assumption implies that the drift term $\tilde{d}(t)$ can be `cancelled out' by choosing the input $\bar{u}(t) = \tilde{B}^\dagger\tilde{d}(t)$. Meanwhile, the parameter $\gamma>0$ in Assumption \ref{asm:stabilizable} measures the difficulty (in terms of a simple $\Ltwo$ cost) of stabilizing the Jacobian Linearizations along system trajectories. 
 
Roughly speaking, our first three conditions ensure that the state costs are in a certain sense `compatible' with the local (first-order) geometry of the control system, meaning that at each point in the optimization space they guide local search algorithms to find an input $\tilde{u}$ which drives the corresponding predictive trajectory $\tilde{x}$ towards the origin. Indeed, by Lemma \ref{lem:comparison}, when Assumptions \ref{asm:cost_smooth}, \ref{asm:noise_structure} and \ref{asm:stabilizable} all hold, at each point $\tilde{u} \in \UT$ the functional $J_T(\cdot;x_0)$ has the same local structure (up to first-order approximations) as an optimal control problem with convex costs and stabilizable time-varying dynamics, namely, $\Jlin(\cdot;x_0,\tilde{u})$. For this local convex approximation it is much clearer to see how the state costs yield state trajectories which decay to the origin. Our second counter-example investigates a situation where there does not exist a coordinate system in which both Assumption \ref{asm:cost_smooth} and Assumption \ref{asm:noise_structure} can be satisfied simultaneous and local search algorithms can produce predictive trajectories which get `stuck' at undesirable equilibria.

Our last technical condition, which is made Assuming \ref{asm:noise_structure} already holds, effectively bounds how costly it is to `cancel out' the affine drift term $\tilde{d}(t)$ along each trajectory:

\begin{assumption}\label{asm:linearizable_bound} \ieeethm
There exists $
L_x,L_u>0$ such that for each system trajectory $(\tilde{x}(\cdot),\tilde{u}(\cdot))$ defined on $\sp{0,T}$ we have   $\| \tilde{B}^\dagger(t)\tilde{d}(t)\| \leq L_x\| \tilde{x}(t)\| + L_u\| \tilde{u}(t)\|$ for each $t \in \sp{0,T}$. Moreover, these constants satisfy
\begin{equation}\label{eq:growth}
L_u^2 \le \tfrac{\alpha_R}{8\beta_R} \ \ \ \text{ and } \ \ \ L_x^2 \le \tfrac{\alpha_Q}{8\beta_R}
\end{equation}
\end{assumption}

In particular, Assumption \ref{asm:linearizable_bound} sates that the cost of rejecting $\tilde{d}(t)$ can only grow linearly with $\tilde{x}(t)$ and $\tilde{u}(t)$. The constant $L_x$ can be made arbitrarily large by re-scaling the relative magnitudes of the state and input costs (so that $\alpha_Q \gg \beta_R$). However, since $\alpha_R \leq \beta_R$, the condition implies that $L_u$ can be at most $\frac{1}{\sqrt{8}}$, which effectively limits how nonlinear the control system is with respect to the input.\footnote{ We remark that the factors of $\frac{1}{8}$ in \eqref{eq:growth} can be replaced by any constant in the interval $(0,1)$ and the proof of Theorem \ref{thm:plan} will go through with minor modifications. However fixing a specific constant simplifies the statement of the main result.} As our first counter example demonstrates, when Assumption \ref{asm:linearizable_bound} is violated local search algorithm may again get `stuck' at undesirable stationary points. The intuition for this failure mode is that even when Assumptions \ref{asm:noise_structure} and \ref{asm:stabilizable} are satisfied if $\tilde{d}(t)$ grows too quickly it may appear `too costly' (from the perspective of optimization algorithms which only have access to first-order information) to reject $\tilde{d}(t)$ and drive the system to the origin. 

Under these assumptions we obtain our main result:
\begin{theorem}\label{thm:plan} \ieeethm Suppose \Cref{ass:v_field_regularity,ass:eq_point,asm:cost_smooth,asm:stabilizable,asm:noise_structure,asm:linearizable_bound} hold. Then for each $T\geq 0$, $x_0 \in \R^n$ and every $\tilde{u} \in \UT$  which is an $\epsilon$-$\FOS$ of $J_T(\cdot;x_0)$ the following hold. If $\alpha_V > 0$ , then $\forall s \in [0,T]$,
\begin{align*}
 ~~\|\tilde{x}(s)\|^2 \le  \calC_0 \cdot ( \calC_1 e^{- \frac{s}{ \calC_1}} \cdot \|x_0\|^2 + \calC_2 \epsilon^2),
\end{align*}
where $\calC_{0} =  6L_F+\tfrac{2L_F\alpha_Q }{ \alpha_R} + \tfrac{2\alpha_Q}{ \alpha_V}$, $\calC_{1}  = 4\gamma\max\{\beta_V,\beta_R,\beta_Q\}$ and $\calC_{2}  =  \tfrac{1}{2\alpha_R}(1+8\beta_R\max\{\tfrac{L_x^2}{\alpha_Q},\tfrac{ L_u^2}{\alpha_R}\})$. More generally, for any $\alpha_V \ge 0$, it holds that for all $s \in \sp{0,T}$ and $\delta > 0$,
\begin{equation*}
 \|\tilde{x}(s)\|^2 \le  \calC_0^\delta\cdot (\calC_1^\delta e^{- \frac{s}{ \calC_1}} \cdot \|x_0\|^2 + \calC_2 \epsilon^2)
\end{equation*}
where $\calC_{0}^\delta :=  6L_F+\tfrac{2L_F\alpha_Q }{ \alpha_R} + \min\{\tfrac{2}{ \delta},\tfrac{2\alpha_Q }{ \alpha_R}\}$,  $\calC_1^\delta := e^{\frac{\delta}{\calC_1}}\calC_1$. 
\end{theorem}

Note that when $\epsilon =0$ taking the square root of both sides of either bound in the statement of the theorem demonstrates that the stationary point is exponentially decaying. We emphasize that the rate of decay is uniform across all stationary points corresponding to different initial conditions $x_0\in \R^n$ and prediction horizons $T>0$. This uniformity is essential for our RHC stability results in Section \ref{sec:FORHC}. There we stipulate how long the prediction horizon $T>0$ and how small the optimality parameter $\epsilon>0$ must be each time a planning problem is solved to ensure stability.

\subsection{Proof of \Cref{thm:plan}}\label{subsec:main_proof}
Let $(\tilde{x}(\cdot), \tilde{u}(\cdot))$ be as in the statement of the theorem, and let $\Vtil(s):= J_{T-s}(\tilde{u}_{[s, T]} ;\tilde{x}(s))$ for $s \in \sp{0,T}$. Then the following bounds hold (under the assumptions of \Cref{thm:plan}):
\begin{lemma}\label{lem:underbound}
If $\alpha_V>0$ then for each $0 \le s' \le s \le T$ we have $\|\tilde{x}(s)\|^2 \le \tfrac{1}{\alpha_Q} \calC_0 \cdot \Vtil(s')$. Alternatively, if $\alpha_V \geq 0$ then for each $0 \leq \delta \leq T$, $s' \leq T-\delta$ and $s' \leq s \leq T$ we have $\|\tilde{x}(s)\|^2 \le \tfrac{1}{\alpha_Q} \calC_0^\delta \cdot \Vtil(s')$.
\end{lemma}
\begin{lemma}\label{lem:overbound} If $\util(t)$ is an $\epsilon$-$\FOS$ of $J_T(\cdot,x_0)$, then for each $s \in \sp{0,T}$ we have $\Vtil(s) \le  \alpha_Q( \calC_1 \|\xtil(s)\|^2 + \tfrac{\calC_2}{2} \epsilon^2)$.  
\end{lemma}
\iftoggle{fullversion}
{
Proofs of both lemmas above are given in the appendix.
}
{

The proof of \Cref{lem:overbound} can be found in the Appendix, while the proof of \Cref{lem:underbound} can be found in \cite{techreport}. However, we note that this under bound can be produced by under bounding the cost-to-go for $J_T(\cdot,x_0)$, and that similar under bounds have appeared in the literature \cite{jadbabaie2001unconstrained}.
}

By the Fundamental Theorem of Calculus,
\begin{align*}
-\tfrac{\rmd}{\rmd s} \Vtil(s) &= Q(\xtil(s)) + R(\util(s)) \ge \alpha_Q \|\xtil(s)\|^2 \\
&\ge  \tfrac{1}{ \calC_1}\Vtil(s) -  \tfrac{\alpha_Q\calC_2 \epsilon^2}{2 \calC_1},
\end{align*}
where the last line uses \Cref{lem:overbound}. Integrating the bound and again invoking \Cref{lem:overbound}, 
\begin{align}\label{eq:exp_value_bound}
\Vtil(s) &\le \exp( - \tfrac{s}{ \calC_1}) \Vtil(0) + \tfrac{\alpha_Q\calC_2 \epsilon^2}{2 \calC_1}\cdcstyle \int_{t=0}^s \exp( - \tfrac{t}{ \calC_1})\rmd t\\ \nonumber
&\le \exp( - \tfrac{s}{ \calC_1}) \Vtil(0) + \tfrac{\alpha_Q}{2}\calC_2 \epsilon^2\\  \nonumber
&\le \alpha_Q \cdot ( \calC_1 e^{- \frac{s}{ \calC_1}} \cdot \|\xtil(0)\|^2 + \calC_2 \epsilon^2).
\end{align}
Finally, in the case where $\alpha_V >0$, \Cref{lem:underbound} lets us convert the above bound to one on $\|\xtil(s)\|^2$, replacing $\alpha_Q$ with $\calC_0$, as desired. In the case where $\alpha_V=0$, for each $ \delta \in \sp{0,T}$ application of \Cref{lem:underbound} yields the desired result for each $s \in \sp{0,T-\delta}$. For each $s \in \sp{T-\delta, T}$ \Cref{lem:underbound} yields
\begin{align*}
\|\tilde{x}(s) \|^2_2 &\leq  \calC^\delta_0 \cdot ( \calC_1 e^{- \frac{T-\delta}{ \calC_1}} \cdot \|\xtil(0)\|^2 + \calC_2 \epsilon^2)\\
& \leq  \calC^\delta_0 \cdot (e^{\frac{\delta}{\calC_1}} \calC_1 e^{- \frac{s}{ \calC_1}} \cdot \|\xtil(0)\|^2 + \calC_2 \epsilon^2).
\end{align*}

\subsection{Connection to Feedback Linearization}\label{subsec:fbl}
We begin by applying our sufficient conditions to feedback linearizable systems, perhaps the most widely studied and well-characterized class of systems in the nonlinear geometric control literature \cite[Chapter 9]{sastry2013nonlinear}. Roughly speaking, a system is feedback linearizable if it can be transformed into a linear system using state feedback and a coordinate transformation. 
Formally, we say that \eqref{eq:dynamics} is feedback linearizable if it is both control-affine, namely, of the form
\begin{equation*}
 F(x,u)= f(x) +g(x)u,
\end{equation*} 
where $f \colon \R^n \to \R^n$ and $g \colon \R^n \to \R^{n\times m}$, and if  there exists a change of coordinates $\xi = \Phi(x)$, where $\Phi \colon \R^n \to \R^n$ is a diffeomorphism, such that in the new coordinates the dynamics of the system are of the form
\begin{equation}\label{eq:linearized_form}
\dot{\xi} = 
\hat{A} \xi + \hat{B}[\hat{f}(\xi) + \hat{g}(\xi)u]
:=\hat{F}(\xi,u),
\end{equation}
where $\hat{A} \in \R^{n \times n}$ and $\hat{B} \in \R^{n\times m}$ define a controllable pair $(\hat{A},\hat{B})$, $\hat{f} \colon \R^n \to \R^m$ and $\hat{g}\colon \R^n \to \R^{n \times m}$ is such that $\hat{g}(\xi)$ is invertible for each $\xi \in \R^n$. We will let $\hat{g}_i(\xi)$ denote the $i$-th column of $\hat{g}(\xi)$. We emphasize that this \emph{global} transformation is distinct from the local Jacobian linearizations employed earlier. In this case the application of the feedback rule $u(\xi,v) = \hat{g}^{-1}(\xi)[-\hat{f}(\xi) + v]$, where $v \in \R^m$ is a new `virtual' input, results in $\dot{\xi} = \hat{A}\xi + \hat{B}v$. In essence, feedback linearization reveals a linear structure underlying the global geometry of the system. Clearly $\hat{F}$ satisfies \Cref{asm:noise_structure}, and the following proposition provides sufficient conditions for \Cref{asm:stabilizable,asm:linearizable_bound} to hold in the new coordinates:
\begin{proposition}\label{prop:fbl1} \ieeethm
Suppose that \eqref{eq:dynamics} is feedback linearizable, and let $\hat{f}$, $\hat{g}$, $\hat{A}$ and $\hat{B}$ be as defined above. Assume that $i)$ there exists $\hat{L} >0$ such that $\| \frac{d}{d\xi}\hat{f}(\xi)\| <\hat{L}$ for each $\xi \in \R^n$ and $ii)$ $\hat{g}(\cdot)$ is constant on $\R^n$. Then there exists $\gamma>0$ such that along each trajectory $(\tilde{\xi}(\cdot), \tilde{u}(\cdot))$ of $\hat{F}$ the associated Jacobian linearization $(\tilde{A}(\cdot), \tilde{B}(\cdot))$ is $\gamma$-stabilizable. Furthermore the drift term $\tilde{d}(t) = \hat{f}(\xi(t)) - \tilde{A}(t)\tilde{\xi}(t)$ satisfies $\|\tilde{B}^\dagger(t) \tilde{d}(t) \|_2 \leq  2\hat{L}\| \tilde{\xi}(t)\|_2$.
\end{proposition}

\begin{remark}\ieeethm
Suppose that the representations of the state running and terminal costs, $\hat{Q} :=Q \circ \Phi^{-1}$ and $\hat{V} :=V\circ \Phi^{-1}$, are convex in the linearizing coordinates and satisfy pointwise bounds as in \Cref{asm:cost_smooth}. Further assume that the assumptions made of $\hat{F}$ in Proposition \ref{prop:fbl1} hold. Then the conclusions of \Cref{thm:plan} can be applied to the representation of $J_T(\cdot,\xi_0)$ in the linearizing coordinates by rescaling $\hat{Q}$ and $R$ appropriately. 
\end{remark}

Thus, the global linearizing coordinates provide a useful tool for verifying the sufficient conditions in \Cref{thm:plan}. Moreover, as we illustrate with our counter-examples, they also provide insight into what goes wrong in cases where the state costs do not lead to stabilizing behavior.

While we illustrate this point further with our examples, let us briefly remark on the necessity of the conditions in Proposition \ref{prop:fbl1} (for our analysis). First, note that along a given solution $(\tilde{\xi}(\cdot), \tilde{u}(\cdot))$ we have $\tilde{A}(t) = \hat{A} + \hat{B}[\frac{d}{d\xi}\hat{f}(\tilde{\xi}(t)) + \sum_{i=1}^{m}\frac{d}{d\xi}g_i(\tilde{\xi}(t))\tilde{u}_i(t)]$ and $\tilde{B}(t) = \hat{B}\hat{g}(\tilde{\xi}(t))$ and $\tilde{d}(t) = \tilde{A}(t)\tilde{\xi}(t) -\frac{d}{d\xi}\hat{f}(\tilde{\xi}(t))\tilde{\xi}$. When condition $ii)$ is violated, even when $\frac{d}{d\xi}\hat{g}(\xi)$ can be bounded globally, the linear growth of $\tilde{A}(t)$ with respect to $\tilde{u}(t)$ may make the pair $(\tilde{A}(\cdot),\tilde{B}(\cdot))$ more difficult to stabilize (in the sense of Assumption \ref{asm:stabilizable}) for large values of the input. Moreover, in this case $\tilde{d}(t)$ will have quadratic cross-terms in $\tilde{\xi}(t)$ and $\tilde{u}(t)$, which may violate the growth conditions in Assumption \ref{asm:linearizable_bound}. Similar issues arise when $\frac{d}{d\xi} \hat{f}(\xi)$ is not bounded globally. This occurs, for example, in Lagrangian mechanical systems wherein the Coriolis terms display quadratic growth in the generalized velocities of the system. The core challenge in each of these cases, from the perspective of our analysis, is that without making additional structural Assumptions beyond those in Proposition \ref{prop:fbl1} it is difficult to rule out cases where the time-varying approximation to the dynamics along some trajectory of the system is arbitrarily difficult to stabilize.

\subsection{Counterexamples}\label{subsec:examples}
\fullvtext{We present three counterexamples and one positive example to illustrate the necessity of our conditions and generality of our results. The first counterexample demonstrates the necessity of small control costs. The second demonstrates the necessity of the matching condition between the drift $\tilde{d}(t)$ and range of the linearized $B$-matrix $\tilde{B}(t)$. The third counterexample constructs an example where our conditions hold, and thus first-order stationary points are stabilizing, but where there exists a $\FOS$ which is \emph{not a global optimum}. This reveals that our conditions are more generaly than popular forms of ``hidden convexity'' such as quasi-convexity. Finally, our fourth example describes a system which is not input affine, but for which the regularity conditions outlined above still hold.}

\begin{figure}[t]
      \centering
      \includegraphics[width = .8\columnwidth]{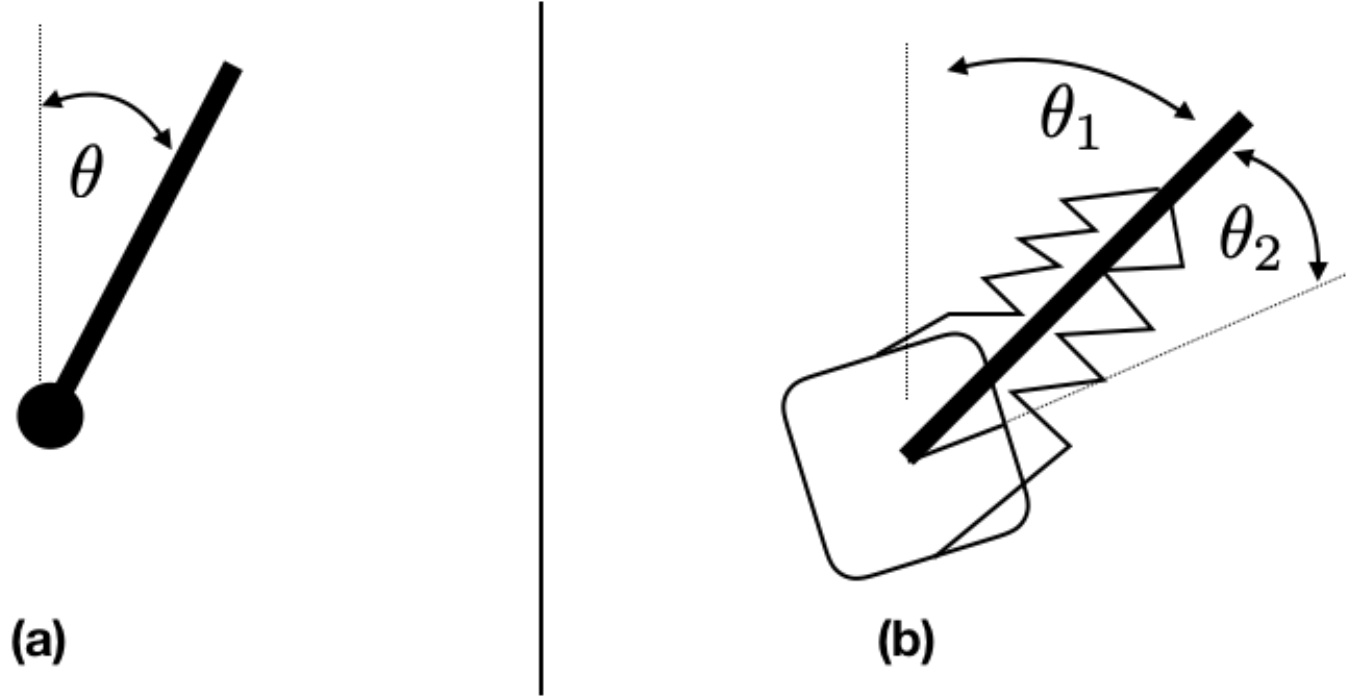}
      \caption{(a) Schematic for the simple inverted pendulum (b) schematic for the inverted pendulum with a flexible joint.}
      \label{fig:examples}
\end{figure}

\cdcparsssec{Relative Weighting of State and Input Costs}
We first consider the simple inverted pendulum in Figure \ref{fig:examples}. The states are $(x_1,x_2)= (\theta,\dot{\theta})$, where $\theta$ is the angle of the arm from vertical. The dynamics are governed by
\begin{align*}
   \begin{bmatrix}
    \dot{x}_1 \\
    \dot{x}_2
    \end{bmatrix} = \begin{bmatrix}
    x_2 \\ 
    k \sin(x_2) + u
    \end{bmatrix},
\end{align*}    
where $k =g\ell$ with $g>0$ the gravitational constant and $\ell>0$ the length of the arm.
The cost to be minimized is
\begin{equation*}
    J_T(\cdot,x_0) = \cdcstyle \int_{0}^T\| \tilde{x}(t)\|_{Q}^2 + r\| \tilde{u}(t)\|_2^2 \dt + \|\tilde{x}(T) \|_{Q_V}^2,
\end{equation*}
where the scalar parameter $r>0$ is used to control the relative weighting of the input and state costs and 
\begin{equation*}
    Q_V = \begin{bmatrix}
    1/4 & 1/(5\sqrt{2})\\
    1/(5\sqrt{2}) & 1
    \end{bmatrix}  \ \ \  Q = \begin{bmatrix}
    1 & -1/4\\
    -1/4 & 1
    \end{bmatrix}
\end{equation*}
Note that the dynamics are linearizable, as they are already in the form \eqref{eq:linearized_form}. Moreover, applying Proposition \ref{prop:fbl1}, one can show that that Assumptions \ref{asm:cost_smooth} through \ref{asm:linearizable_bound} are satisfied with parameters $\alpha_Q = \beta_Q = \frac{3}{4}$, $\alpha_R = \beta_R =r$, $L_x = k$ and $L_u =0$. Considering $k = 10$, we find that if if $r < \frac{3}{4k^2} = \frac{3}{400}$ the sufficient conditions for exponential stability of Theorem \ref{thm:plan} are satisfied. However,  if $r$ is not small enough then there may exists undesirable first order stationary points of the cost functional. Specifically, consider the initial condition $x_0 = (\frac{3 \pi}{4}, 0)^T$ and the control signal $\tilde{u}(\cdot) \equiv -k\sin(\frac{3 \pi}{4}) = -\frac{10}{\sqrt{2}}$, which generates the trajectory $\tilde{x}(\cdot) \equiv x_0$. The costate along this arc is $p(\cdot) \equiv (\frac{3 \pi}{16},\frac{3 \pi}{20\sqrt{2}})$ and the gradients of the objective is given by $\nabla J_T(\tilde{u};x_0)(t) = p(t) + ru(t)$. Thus, we see that if we choose $r = \frac{3 \pi}{200}$ then we will have $\nabla J_T(\tilde{u};x_0)(\cdot) \equiv 0$, which demonstrates that $\tilde{u}$ is a stationary point of the cost function. Thus, while RHC stability results which rely on global optimality of each planning problem predict stabilizing behavior for sufficiently large $T>0$ \cite{jadbabaie2001stability}, this example demonstrates that algorithms which only find first order-stationary points may be `too myopic' to guarantee stability unless the input cost is small enough.

\fullvtext{
It is interesting to note that it is impossible to find a stationary pair $(\hat{x}(\cdot),\hat{u}(\cdot))$ of $J_T(\cdot,x_0)$ with the property that $\hat{x}(\cdot) \equiv x_0$ if we instead have $x_0 \in (-\frac{\pi}{2}, \frac{\pi}{2})$. Indeed, if we pick $\tilde{u}(\cdot) = -\sin(x_0)$ so that $\hat{x}(\cdot) \equiv x_0$, then in this case the costate will satisfy the differential equation $-\dot{p}(t) = \hat{A}(t)p(t) + 2q\hat{x}(t)$ where $\hat{A}(\cdot) \equiv  \cos(x_0) > 0$. Thus, in this case the costate cannot be a constant function of time, which means that $\hat{u}(\cdot)$ cannot be a stationary point of $J_T(\cdot,x_0)$. 

More broadly, consider a general 1-dimensional system $\dot{x} =F(x,u)$ which satisfies \Cref{asm:stabilizable,asm:noise_structure}. One can verify that and trajectory $(\xcheck(\cdot), \ucheck(\cdot))$ of $F$ such that $\xcheck(\cdot) \equiv x_0$ and $\ucheck(\cdot) \equiv \ucheck_0$ (for some $x_0, \ucheck \in \R$) can only be a stationary pair of $J_T(\cdot,x_0)$ if $\frac{\partial }{\partial x}F(x_0, \bar{u})<0$ and $\check{d}(\cdot) \equiv \check{d}_0 = F(x_0, \ucheck_0) - \check{A}x_0 -\check{B}u_0$ is such that $sign(d_0) = sign(x_0)$, where $(\check{A}_0,\check{B}_0)$ are the Jacobian linearization of $F$ at $(x_0,\ucheck_0)$. This highlights the need to further study how the geometry of the Jacobian linearizations and drift terms promote or inhibit stabilizing behavior.}

\newcommand{\calD}{\mathcal{D}}
\cdcparsssec{Structure of Local Drift Term}
Consider the flexible link manipulator depicted in Figure \ref{fig:examples}. The state is $(x_1,x_2,x_3,x_4)^T = (\theta_1,\dot{\theta}_1,\theta_2,\dot{\theta}_2)^T$, where $\theta_1$ is the angle of the arm from vertical and $\theta_2$ is the difference between the angle of the arm and the internal angle of the motor. The dynamics are
\begin{equation*}
     \begin{bmatrix}
    \dot{x}_1 \\
    \dot{x}_2 \\
    \dot{x}_3 \\
    \dot{x}_4
    \end{bmatrix} = \begin{bmatrix}
    x_2 \\
    -K_1 \sin(x_1) + K_2 x_3 \\
    x_4\\
    -K_1 \sin(x_1) - (K_2 + K_3)x_3 - u
    \end{bmatrix},
\end{equation*}
where $K_1 = g\ell$, $K_2 = \frac{k}{M}$ and $K_3 =\frac{k}{I}$, where $g$ is the gravitational constant, $\ell$ is the length of the arm, $k = 10$ is the spring coefficient, $M$ the mass of the arm and $I$ the internal inertia of the motor. For concreteness, we will assume that these physical parameters are such that $K_1 = 20$ and $K_2 = K_3 = 1$.  We apply cost functionals of the form
\begin{align*}
    J_T(\cdot ; x_0) =\cdcstyle \int_{0}^T \| \tilde{x}(t)\|_{Q} + r\| \tilde{u}(t)\|_2^2 \dt + \|\tilde{x}(T)\|_{Q_V},\\
    Q:= \begin{bmatrix}
    1 &\frac{1}{4} & 0 & 0 \\
    \frac{1}{4} & 1 & 0 & 0\\
    0 & 0 & 1 & 0 \\
    0 & 0 & 0 & 1
    \end{bmatrix} \ \ \ Q_V:= \begin{bmatrix}
    5 & -1 & -\frac{1}{4} & 0 \\
    -1 & 5 & 0 & 0\\
    -\frac{1}{4} & 0 & 5 & 0 \\
    0 & 0 & 0 & 5
    \end{bmatrix}.
\end{align*}
While the running and terminal costs satisfy Assumption \ref{asm:cost_smooth} in the $x$ coordinates, Assumption \ref{asm:noise_structure} is violated for this parameterization of the control system. Consider the initial condition $x_0 = (\pi, 0,0,0)^T$. Note that the input $\tilde{u}(\cdot) \equiv 0$ generates the trajectory $\tilde{x}(\cdot) \equiv x_0$. The reader may verify that the costate along this trajectory is $p(t) \equiv (\pi, \frac{1}{4}\pi,0,0)$ and the gradient at this point in the optimization space is defined by $J_T(\tilde{u};x_0)(\cdot)\equiv 0$. Note that this is true for every choice of prediction horizon $T>0$ and choice of the scaling parameter $r>0$. Thus, regardless of the prediction horizon, algorithms which find first-order stationary points may get stuck at this undesirable equilibrium. 

We can also see how the proposed cost function fails to guide local search algorithms to stabilizing solutions by studying its structure in a set of linearizing coordinates. Indeed, consider coordinates defined by $\xi = \Phi(x)$ where  
\begin{equation*}
   \Phi(x) = (x_1, x_2,-K_1\sin(x_1) + K_2 x_3, -K_1\cos(x_1)x_2 + K_2 x_4).
\end{equation*}
These coordinates are obtained by input-output linearizing the dynamics with the output $y = x_1$ (see \cite[Chapter 9]{sastry2013nonlinear}), and in the new coordinates the system can be shown to be of the form \eqref{eq:linearized_form}. Thus, in the new coordinates Assumption \ref{asm:noise_structure} is satisfied, however, Assumption \ref{asm:cost_smooth} is not satisfied. Indeed, due to the nonlinearities in $\Phi$, the reader may verify that the maps $z \to \|\Phi^{-1}(z) \|_{Q}^2$ and $z \to \|\Phi^{-1}(z)\|_{Q_V}$ are not convex. Thus, in these coordinates, the local structure of the state costs attracts trajectories towards the undesirable equilibrium point we identified above.

\cdcparsssec{Non-Control-Affine System}
Consider the scalar system $\dot{x} = f(x)+ bu + k\cos(x)\sin(u)$,
where $f \colon \R \to \R$ satisfies $f(0) =0$ and $\|\frac{d}{dx}\tilde{f}(x) \|<L$ for each $x \in \R^n$ and  $b>k>0$. Along a given system trajectory $(\tilde{x}(\cdot),\tilde{u}(\cdot))$ the Jacobian linearization is $(\tilde{A}(\cdot),\tilde{B}(\cdot))$ where $\tilde{A}(t) = \frac{d}{dx}f(\tilde{x}(t)) - k \sin(\tilde{x}(t))\sin(\tilde{u}(t))$ and $\tilde{B}(t) = b + k\cos(\tilde{x}(t))\cos(\tilde{u}(t))$. Note that $|\tilde{A}(t)| < L + k$ and $\tilde{B}(t) > b -k$, and thus it is straightforward to argue that the linearizations along the trajectories satisfy Assumption \ref{asm:stabilizable} for a suitable parameter $\gamma>0$. The disturbance term associated to the chosen trajectory is $\tilde{d}(t) = f(\tilde{x}(t)) +k\cos(\tilde{x}(t))\sin(\tilde{u}(t))  -\tilde{A}(t)\tilde{x}(t) -  k\cos(\tilde{x}(t))\cos(\tilde{u}(t)) \tilde{u}(t)$, which can be seen to satisfy $\|B^\dagger(t) \tilde{d}(t)\| \leq \frac{1}{b-k} ((2L+k)\|\tilde{x}(t)\|  + 2k \| \tilde{u}(t)\|)$. Thus, we see that these dynamics satisfy Assumption \ref{asm:linearizable_bound}. Further, if in addition we have$\frac{2k}{b-k} <\frac{1}{8}$ then Assumption \ref{asm:linearizable_bound} can be satisfied by choosing cost functions which Satisfy Assunption \ref{asm:cost_smooth} and weight the state and input penalties appropriately. While most of the systems we consider in this paper are control-affine, this example demonstrates that \Cref{thm:plan} can be applied to non-control-affine systems so long as the nonlinearity in the input is not `too large'.

\subsection{Stabilizing stationary points need not be global optima. }
In this example, we demonstrate that exponentially stabilizing staionary points need not be global optima. Thus, the problems do not exhibit ``hidden convexity'' or quasi-convexity.

To construct the counterexample, we construct a system which has linear dynamics $\dot{x} = u$ within a region $\calD_1$, but which enjoys a beneficial drift term in a region $\calD_2$, with dynamics $\dot{x} = u + d$. The two regions are interpolated via a bump function. We construct a cost $J_T$ so there is a $\FOS$ which remains in $\calD_1$, executing the optimal control law for the driftless dynamics. We then show that the cost can be strictly improved with a control law which enters $\calD_2$ to take advantage of the drift.

Formally, let $\eta: \R \to \R$ be any bump function with $\eta'$ bounded nonnegative, with $\eta(z) = 0$ for $z \le 1$, and $\eta(z) = 1$ for $z \ge 2$. For convenience, assume the normalization $\int_{1}^2 \eta(z) \rmd z = 1/2$.  Consider the system in $\R^2$ with state $x = (x_1,x_2)$, input $u = (u_1,u_2)$, and, for a parameter $\alpha > 0$ to be manipulated, the dynamics 
\begin{align}
    \dot{x} = u + \eta(x_1)\underbrace{\begin{bmatrix} 0\\
    -\alpha
    \end{bmatrix}}_{d} 
\end{align}
For time horizon $T$ and parameters $q,r,B > 0$, consider an initial point $x_0 = (0,B)$
\begin{align}
    J_T(\tilde{u};x_0) &= \int_{0}^T (q\|\tilde{x}(t)\|_2^2 + r\|\tilde{u}(t)\|_2^2)\dt + \sqrt{qr}\|\tilde{x}(T)\|_2^2 \label{eq:J_T} \\
    &\text{s.t.} \quad \dot{x} = u + \eta(x_1) \cdot d \nonumber.
\end{align}
Much like the previous example, one can verify that the above system satisfies the requisite stability conditions to ensure that all first-order stationary points are exponentially stabilizing. We now show that if $B$ satisfies
\begin{align}
\frac{B}{18} \ge \max\left\{1 + \frac{1}{\alpha \sqrt{rq}},\sqrt{\frac{r}{q}}\right\} \label{eq:Bchoice},
\end{align}
and if $T$ is sufficiently large, then there exists a local optimum of $J_T$ which is not global. 

To show this, consider the region $\calD_1 = \{(x_1,x_2) \in \R^2: x_1 < 1\}$. On $\calD_1$, the dynamics are simply $\dot{x} = u$, and thus a first order stationary point which remains in $\calD_1$ is given by the optimal control law from that system, i.e. the optimal law of
\begin{align}
\bar{J}_T(\tilde{u};x_0) =  \int_{0}^T (q\|\tilde{x}(t)\|_2^2 + r\|\tilde{u}(t)\|_2^2)\dt + \sqrt{qr}\|\tilde{x}(T)\|_2^2  \quad \text{s.t.} \quad \dot{x} = u. \label{eq:barJ_example}
\end{align}
To classify the optimal control law, observe that the optimal law for $\lim_{T \to \infty}  J_T(\tilde{u};x_0)$ can be obtained by solving the Continuous Algebraic Riccati Equation with matrices $A = 0$, $B = I$, yielding a cost-to-go matrix $P = \sqrt{qr}I$, and optimal law $\tilde{u}^\star = -\sqrt{q/r}\cdot\tilde{x}^\star$. Thus, the terminal cost in $J_T$ is equal to the infinite-horizon cost-to-go function $\tilde{x}(T)^\top P \tilde{x}(T)$, $\tilde{u}^\star(t) = -\sqrt{q/r}\cdot\tilde{x}^\star(t)$ is an optimal law for $\bar{J}_T$ as in \eqref{eq:barJ_example} with value $\bar{J}^\star_T(x_0) = \inf_{\tilde{u}}\bar{J}_T(\tilde{u};x_0) = \sqrt{rq}B^2$. Note that the law $\tilde{u}^\star(t) = -\sqrt{q/r}\cdot\tilde{x}^\star(t)$, when executed starting at $x_0 = (0,B)$, remains on the Y-axis in $\R^2$, and hence remains in $\calD_2$. Thus, this law is also a local optimal of $J_T(\cdot;x_0)$ with cost $\sqrt{rq}B^2$.

To show that $\tilde{u}^\star(t) = -\sqrt{q/r}\cdot\tilde{x}^\star(t)$ is not a global minimizer of $J_T(\cdot;x_0)$, let us construct a control $\hatu(t)$ with $J_T(\hatu;x_0) < J_T(\util^\star;x_0)$. Let $t_0 = 1/(B \sqrt{q/r})$,  $t_1 = \alpha(1 - t_0)$ to be selected, consider the control $\hatu(t)$ 
\begin{align}
\hatu{u}(t) = \begin{cases} (\frac{2}{t_0},0) & t \in [0,t_0]\\
(0,0)) & t \in [t_0,t_0+t_1]\\
(-\frac{2}{t_0},0),0) & t \in [t_0+t_1,2t_0 + t_1]\\
(0,0) & t \ge t_1 + 2t_0
\end{cases}
\end{align}
Let us characterize the trajectory of $\hat{x}(t)$:
\begin{itemize}
\item For $t \in [0,t_0]$, $\hat{x}_1(t)  = 2t/t_0$, and $\hat{x}_2(t)$ is monotonically decreasing (since $\eta(z)$ is nonnegative) and gives $\hat{x}_2(t_0) = B - \int_{0}^{t_0} \eta(\hat{x}_1(t))\dt = B - \int_0^{t_0}\eta(2t/t_0)\dt = B - 2t_0 \int_0^2 \eta(z)\rmd z = B - 2t_0 \cdot \frac{1}{2} = B - t_0$. The total cost incurred during this time is at most 
\begin{align}
t_0(q2^2 + qB^2 + r\frac{4}{t_0^2})  = \frac{4q}{B\sqrt{q/r}} + \frac{qB^2}{B\sqrt{q/r}} + 4r B\sqrt{q/r} = \sqrt{rq}(4/B + 5B).
\end{align}
\item For $t \in [t_0,t_1+t_0]$, $\hat{x}_1(t) = 2$, and $\ddt \hat{x}_2(t) = -\alpha $; $\hat{x}_2(t_0 + s) = B - t_0 - \alpha s$. For our choice of $t_1 = \frac{1}{\alpha}(B - 2t_0)$, we have $\hat{x}_2(t_1+t_0) = t_0$. The total cost incurred during the time is at most $t_1 (4 + B^2) \le (4/B + B)/\alpha \le 4(1/B + B)/\alpha$.
\item One can verify that, for $t \in [t_1+t_0,t_1+2t_0]$, $\hat{x}_1(t) \in [0,2]$ and  $\hat{x}_2(t_1+2t_0) \in [0,t_0]$. Moreover, $\hat{x}(t_1+2t_0) = (0,0)$. The total cost incurred is at most $q(4 + t_0^2)t_0 + 4r/ t_0$. The choice of $B$ as in \eqref{eq:Bchoice} ensures $t_0 \le 1$, so that we can upper bound the incurred cost by $5qt_0 + 4r/t_0$, which is at most
\begin{align}
\sqrt{rq}(5/B + 4B).
\end{align}
\item For the remainder of the trajectory, we remain at $(0,0)$ and incur no cost.
\end{itemize}
Since \Cref{eq:Bchoice}  ensures $B \ge 1$, the total cost is at most
\begin{align}
J_T(\hat{u};x_0) \le (4/B + B)/\alpha + 9\sqrt{rq}(1/B + B) \le 18B(\sqrt{rq}+1/\alpha).
\end{align}
In sum, we have $J_T(\hat{u};x_0) < J_T(\tilde{u};x_0)$ as soon as $18B(\sqrt{rq}+1/\alpha) \le B^2\sqrt{rq}$; that is, as soon as \eqref{eq:Bchoice} holds.


\section{First-Order Stability Guarantees for Receding Horizon Control \label{sec:rhc}}
\newcommand{\CRHC}{\mathtt{CRHC}}
\newcommand{\RHC}{\mathtt{RHC}}
\newcommand{\IRHC}{\mathtt{I}\text{-}\mathtt{RHC}}
\newcommand{\FORHC}{\mathtt{FO}\text{-}\mathtt{RHC}}
\newcommand{\PCshort}{\mathcal{PC}}
\newcommand{\plan}{\mathtt{plan}}
\newcommand{\uplan}{u^{\plan}}
\newcommand{\uwarm}{\bar{u}}
\newcommand{\boldx}{\mathbf{x}}
\newcommand{\boldu}{\mathbf{u}}
\newcommand{\hatxnot}{\hat{x}_0}

\newcommand{\shrink}{\scal}

In \emph{receding horizon control} ($\RHC$) or \emph{model predictive control}, a planner solves $ \inf_{\tilde{u}(\cdot)}J_T(\tilde{u},x(t))$, where $x(t)$ is the current state of the real world system. The planner then applies the resulting open loop predictive control until a new state measurement is received and the process can be repeated. As discussed above, most formal stability guarantees require that an (approximate) globally optimal solution is found for each (generally nonconvex) planning problem. Applying \Cref{thm:plan}, we provide the first stability guarantees for a formal model of nonlinear RHC which only requires planning to approximate stationary points.


\subsection{First-Order Receding Horizon Control \label{sec:FORHC}}

Many practical implementations of RHC use a technique known as \emph{warm starting}, where the predictive control returned during each optimization phase is used to construct the `initial guess' for the subsequent planning problem. This approach has proven highly effective for systems which require rapid re-planning to maintain stability \cite{bock1999efficient}. 

To model this approach, we define the \emph{first-order receding horizon control} strategy, denoted $\FORHC$, as follows. First a prediction horizon $T> 0$ and a replanning interval $\delta \in (0,T]$ are chosen and a sequence of replanning times $t_k = k\delta$ for $k \in \mathbb{N}$ are defined. Next, the process takes in an initial condition of the physical system $x_0 \in \R^n$ and a warm-start  control $\uwarm_0 \in \UT$ specified by the user.  We let  $(\boldx(\cdot,x_0,\uwarm_0), \boldu(\cdot,x_0,\uwarm_0))$ denote the resulting trajectory produced by the control scheme described below. 

At each $t_k$ for $k \in \N$ a warm-start routine generates an initial guess $\bar{u}_k(\cdot) = \uwarm_k(\cdot; x_0,\uwarm_0)\in \UT$ for the problem $J_T(\cdot;\boldx(\cdot,x_0,\uwarm_0))$; a simple choice for such a routine is presented momentarily. The local search method then optimizes the problem using the chosen initial guess, and produces the predictive control $\tilde{u}_{k}(\cdot) = \tilde{u}(\cdot; x_0, \bar{u}_0)$. Note that both of these quantities depend on both the initial condition of the system and the initial warm-start control specified by the user. The predictive control is constructed via $\tilde{u}_k(\cdot) =\uplan_T( \cdot, \boldx(t_k,x_0,\uwarm_0), \uwarm_k) \in \UT$, where the map $\uplan_T$ is used to model how the chosen search algorithm selects a predictive control given for a given initial condition and warm-start input. Finally, the actual control $ \boldu(t,x_0,\uwarm_0) = \tilde{u}_{k}(t-t_k)$ is applied on the interval $[t
_k,t_{k+1})$, and the process repeats.
\begin{assumption}\label{asm:planner} \ieeethm We assume that, for any $\hatxnot \in \R^n,\hatu \in \UT$, the planned solution $\tilde{u} = \uplan_T(\cdot,\hatxnot,\hatu) \in \UT$ satisfies the following two conditions with parameter $\epsilon_0 > 0$:

\begin{enumerate}
    \item $J_{T}(\tilde{u};\hat{x}_0) \leq J_{T}(\bar{u};\hat{x}_0)$; and,
    \item $\tilde{u}$ is an $\epsilon_0 J_{T}(\tilde{u}; \hat{x}_0 )^{1/2}$-$\FOS$ of $J_T(\cdot~;\hat{x}_0)$.
\end{enumerate}
\end{assumption}

The rationale for the first condition is that many popular trajectory optimization methods are \emph{descent methods}, and therefore only decrease the value of the functional $J_T$. The second condition is reasonable because such methods converge to approximate first-order stationary points, even for nonconvex landscapes \cite{bertsekas1997nonlinear}. The normalization by $J_{T}(\hat{x}_0, \hat{u})^{1/2}$ 
affords geometric stability in \Cref{thm:RHC} by ensuring the optimization terminates close enough to a stationary point for each planning problem as the system trajectory approaches the origin. 


It remains to specify how the warm-starts $\uwarm_k$ are produced for $k \geq1$. We propose selecting $\tilde{u}_k$ with $\delta$-delay, continuing until time $T$, and then applying zero input:
\begin{equation*}
    \bar{u}_{k+1}(t,x_0,\uwarm_0) = \cdcstyle \begin{cases}
    \tilde{u}_k(t+\delta,x_0,\uwarm_0) &  t \in [0,T-\delta] \\
    0 & t \in (T-\delta ,T].
    \end{cases}
\end{equation*}
While more sophisticated warm-starts may be adopted in practice, the above is preferable for the present analysis because (a) it does not require further system knowledge, and (b) is ammenable to transparent stability guarantees.
\newcommand{\calK}{\mathcal{K}}

\subsection{Sufficient Conditions for Exponential Stability of $\FORHC$}

Finally, we apply \Cref{thm:plan} and its assumptions to provide sufficient conditions for the stability of $\FORHC$. In order to obtain exponential convergence, we will require that, given a desired replanning interval $\delta>0$, the prediction horizon $T>0$ is sufficiently large and the optimality parameter $\epsilon_0>0$ in Assumption \ref{asm:planner} is sufficiently small:

\begin{theorem}\label{thm:RHC}\ieeethm
Let the assumptions in \Cref{thm:plan} hold. Further assume that the search algorithm chosen for $\FORHC$ satisfies the  conditions in \Cref{asm:planner}. Then for any prediction horizon $T >0$, replanning interval $\delta \in (0,T]$ and optimality parameter $\epsilon_0<\sqrt{2 \alpha_Q \calC_2}$, and each initial condition for the physical system $x_0 \in \R^n$ and initial warm-start decision variable $\bar{u}_0 \in \UT$ the system trajectory $(\boldx(\cdot,x_0,\uwarm_0), \boldu(\cdot,x_0,\uwarm_0))$ generated by the corresponding $\FORHC$ scheme satisfies
\begin{equation*}
\|\boldx(t_k,x_0,\uwarm_0) \|_2 \leq  \sqrt{M(\delta,T,\epsilon_0)} e^{\eta(\delta,T,\epsilon_0)t_k}\|x_0\|_2,
\end{equation*}
for each $k \in \N$
where we define
\begin{align*}
M(\delta,\epsilon_0) &:={\calC}_{0}^\delta \calC_1\left( 1 -\tfrac{\alpha_Q}{2}\calC_2\epsilon_0^2  \right)^{-1}\\
\eta(\delta,T,\epsilon_0) &:= \tfrac{1}{2\delta}\ln\left( e^{-\delta / \calC_1} + \mathcal{T}(\delta, T)+\mathcal{E}(\delta,\epsilon_0) \right)\\
\mathcal{E}(\delta,\epsilon_0)&:= \tfrac{1}{2}{\calC}_{0}^\delta \,\calC_2 e^{2L_F\delta}((\delta +1)\alpha_Q + \alpha_V)  \epsilon_0^2 \\
\mathcal{T}(\delta, T) &:={\calC}_{0}^\delta \calC_1(\delta \alpha_{Q} + \alpha_V)e^{-\frac{T}{\calC_1} +\delta(\frac{1}{\calC_1} + 2L_F)}.
\end{align*}
\end{theorem}

To interpret the above constants, first note that for a fixed replanning interval $\delta>0$ we have $\lim_{\epsilon_0 \to 0} M(\delta, \epsilon_0) = \bar{\calC}_0^\delta \calC_1$ and $\lim_{ \substack{%
    T \to \infty\\
    \epsilon_0 \to 0}} \eta(\delta, T,\epsilon_0) =\frac{1}{2\calC_1}$. Thus, in the limiting case $\FORHC$ recovers the exponential rate of convergence predicted by \Cref{thm:plan}. Next, note that $\eta(\delta,T,\epsilon_0)$ will only be negative if $e^{-\frac{\delta}{\calC_1}} + \mathcal{T}(\delta ,T, \epsilon_0) + \mathcal{E}(\delta,\epsilon_0) <1$. Thus, for our estimate on the rate of convergence to be exponentially decaying we require that $T$ is (at least) as big as $T  \geq \ln(\bar{\calC}_0^\delta \calC_1(\delta \alpha_Q + \alpha_V) ) + \delta (1 + \calC_12L_F)$.

\section{Future Directions \label{sec:conclusion}}

There are many important directions for future work. First, it should be determined whether the strong assumptions required for \Cref{thm:plan} can be relaxed. Concretely, in the context of control affine systems of the form $F(x,u) = f(x) +g(x)u$, it remains to be determined under what conditions practical RHC methods can stabilize the system when the drift term $\tilde{d}(t)$ does not grow linearly in $\tilde{x}(t)$ and $\tilde{u}(t)$ (which can happen when the growth of $f$ is super-linear or when $g$ is not constant). The primary shortcoming of our proof technique is that we completely `cancel out' $\tilde{d}(t)$ when constructing sub-optimal controls in the proof of Theorem \ref{thm:plan}, and in some cases this cancellation may be too `costly' to obtain a useful upper bound. However, the disturbance term may actually be useful at certain points if it helps drive the system towards the origin, making complete cancellation unnecessary. It is interesting to note that feedback linearizing controllers, which can also cancel out 'useful' nonlinearities, often receive a similar criticism. Thus, we suspect that in general further structural assumptions on the system dynamics could lead to positive results.





\appendix
\iftoggle{fullversion}
{ 

\section{Appendix}
\begin{proof}[\textbf{Proof of \Cref{lem:underbound}}] Under \Cref{ass:v_field_regularity,ass:eq_point}, we have that $\|\ddt \tilde{x}(t)\| = \|F(\tilde{x}(t),\tilde{u}(t))\| \le L_F(\|\tilde{x}(t)\| + \|\tilde{u}(t)\|)$. Hence, $|\ddt \|\tilde{x}(t)\|^2| = |\langle \tilde{x}(t), \ddt \tilde{x}(t) \rangle| \le L_F\|\tilde{x}(t)\|^2 + L_F \|\tilde{u}(t)\|\|\tilde{x}(t)\|$, which  at most $\frac{3L_F}{2}\|\tilde{x}(t)\|^2 + \frac{L_F}{2}\|\tilde{u}(t)\|^2$ by the AM-GM inequality. \Cref{asm:cost_smooth} then implies this is at most $\frac{c_1}{2\alpha_Q} (Q(\tilde{x}(t)) + R(\tilde{u}(t))$, where we define $c_1 := L_F(3+\frac{\alpha_Q }{ \alpha_R})$. Thus, given any times $s_1 \le s_2 \in [0,T]$, we have
\begin{align}
&\left|\|\tilde{x}(s_1)\|^2 - \|\tilde{x}(s_2)\|^2\right| \le  \int_{t=s_1}^{s_2} \left|\tfrac{\rmd}{\rmd t}\|\tilde{x}(t)\|^2\right| \rmd t \nonumber\\
&\quad  \le \tfrac{c_1}{2\alpha_Q}   \int_{t=s_1}^{s_2}(Q(\tilde{x}(t)) + R(\tilde{u}(t))\dt \le \tfrac{c_1}{2\alpha_Q} \Vtil(s_1) \label{eq:underbound_first_bound}.
\end{align}

To conclude, fix a time $s \in [0,T]$. We consider two cases: \textbf{Case 1:} There is a time $\tau \in [s,T]$ such that $\|\tilde{x}(\tau)\|^2 \le \frac{1}{2}\|\tilde{x}(s)\|^2$. Invoking \eqref{eq:underbound_first_bound},
\begin{align*}
&\tfrac{1}{2}\|\tilde{x}(s)\|^2 \ge  \left|\|\tilde{x}(s)\|^2 -  \|\tilde{x}(\tau)\|^2\right| \\
&\quad \ge  \left|\|\tilde{x}(s)\|^2 - \|\tilde{x}(s\tau)\|^2\right| \ge \tfrac{c_1}{2\alpha_Q} \Vtil(s).
\end{align*}
\textbf{Case 2:} There is no such time $\tau$, so $\|\tilde{x}(t)\|^2 \ge \frac{1}{2}\|\tilde{x}(s)\|^2$ for all $t \in [s,T]$. In this case, $\Vtil(s) \ge \alpha_Q \int_{t=s}^T \|\tilde{x}(t)\|^2 \rmd t + \alpha_V\|\tilde{x}(T)^2\| \ge (\alpha_Q(T-s) + \alpha_V)\cdot\|\tilde{x}(s)\|^2/2$. Inverting, and combining both cases,
\begin{align}
\|\tilde{x}(s)\|^2 \le \frac{1}{\alpha_Q}\max\{c_1, \tfrac{2}{T-s + \alpha_V/\alpha_Q}\}\Vtil(s)  \label{eq:intermediate}
\end{align} 
Finally, for any $s' \in [0,s]$, arguing as in \eqref{eq:underbound_first_bound}, and applying \eqref{eq:intermediate} and some simplifications (incuding $\Vtil(s') \ge \Vtil(s)$),
\begin{align*}
\|\xtil(s)\|^2 &\le \|\xtil(s')\|^2 + |\int_{t=s'}^{s} (\tfrac{\rmd}{\rmd t}\|\tilde{x}(t)\|^2) \rmd t|\\
&\le \|\xtil(s')\|^2 + \frac{c_1}{2\alpha_Q}\Vtil(s')\\
&\le \frac{1}{\alpha_Q}\underbrace{(2c_1 + 2\tfrac{1}{T-(s') + \frac{\alpha_V}{\alpha_Q}})}_{:=\calC_0(s')}\Vtil(s'),
\end{align*}
which can be specialized to the desired cases. 
\end{proof}

\begin{proof}[\textbf{Proof of \Cref{lem:overbound}}]
In view of \Cref{lem:comparison}, to obtain a bound on $\Vtil(s)$ it suffices to bound  $\Vjacst(s) :=\inf_{\ubar_{[s,T]}} \Jlin_{T-s}(\cdot, ~ \xtil(s), \util_{[s,T]})$. Moreover, we can bound $\Vjacst(s)$ by bounding $\Jlin_{T-s}(\ubar_{[s,T]};~ \xtil(t), \util_{[s,T]})$ for any (possibly suboptimal) control $\ubar_{[s,T]}$; for simplicity, let us drop the ${[s,T]}$-subscript going forward. We select $\ubar(t) = \ubar_1(t) + \ubar_2(t)$, where $\ubar_1(t)$ satisfies $\Btil(t) \ubar_1(t) =- \tilde{d}(t)$, and where $\ubar_2$ witnesses $\gamma$-stabilizability at time $s$ as in \Cref{asm:stabilizable}.

With this choice of $\ubar(t)$ the dynamics of $\xbar(t)$ in $\Jlin_{T-s}$  are $\tfrac{\rmd}{\rmd t} \xbar(t) = \Atil(t)\xbar(t) + \Btil(t)\ubar(t) + \tilde{d}(t) = \Atil(t)\xbar(t) + \Btil \bar{u}_2(t)$, and, writing out $\Jlin_{T-s}$ explicitly, we obtain
\begin{align} 
\Vjacst(s) \le  \int_{t=s}^T (Q(\xbar(t)) + R(\ubar(t))\rmd t + V(\xbar(T)).
\end{align}
By the elementary bound $\|\ubar(t)\|^2 \le 2\|\ubar_1(t)\|^2 + 2\|\bar{u}_2(t)\|^2$, the following holds for constant $c = \max\{\beta_V,2\beta_R,\beta_Q\}$,
\begin{align}\textstyle
&\Vjacst(s)\le   2\beta_R \int_{t=s}^T \|\ubar_1(t)\|^2 \dt \label{eq:P_line_one}\\
&\quad+c \left(\int_{t=s}^T (\|\xbar(t)\|^2 + \|\ubar(t)\|^2)\rmd t + \|\xbar(T)\|^2\right). \label{eq:P_line_two}
\end{align}
To bound \eqref{eq:P_line_two}, we observe that $\ddt \xbar(t) = \Atil(t)\xbar(t) + \Btil(t)\ubar(t) + \tilde{d}(t) = \Atil(t)\xbar(t) + \Btil(t)\bar{u}_2(t)$, which corresponds to the $\hat{x}(t)$ dynamics in the definion of $\gamma$-stabilizability; thus, \eqref{eq:P_line_two} is at most $c \cdot \gamma \|\xtil(s)\|^2$. 

To bound \eqref{eq:P_line_one}, we use \eqref{asm:linearizable_bound} to bound $ \|\ubar_1(t)\|^2 \le  2L_x^2\|\xtil(t)\|^2 + 2 L_u^2\|\util(t)\|^2 \le 2c'(Q(\xtil(t)) +  R(\util(t))/\beta_R$, where $c' =  \beta_R\max\{\frac{L_x^2}{\alpha_Q},\frac{ L_u^2}{\alpha_R})$. Hence,  in view of \Cref{lem:comparison}, and the principle of optimality:
\begin{align*}
& 2 \beta_R\textstyle\int_{t=s}^T \|\bar{u}_1(t)\|^2 \dt \le  4c'\beta_R \textstyle\int_{t=s}^T (Q(\xtil(t)) +  R(\util(t))\dt \\
&\quad \le  4c' \beta_R \Vtil(s)  \le 4c'\beta_R(\Vjacst(s)  + \tfrac{\epsilon^2}{2\alpha_R})
\end{align*}
Putting the bounds together and rearranging:
\begin{align*}
(1-   4 \beta_R c' )\Vjacst(s) \le c \cdot \gamma \|\xtil(s)\|^2 + \tfrac{2 c'\beta_R}{ \alpha_R} \epsilon^2
\end{align*}
Under \Cref{asm:linearizable_bound}, we have $4\beta_R c'\le 1/2$, so that
\begin{align}
\Vjacst(s) \le 2c\cdot \gamma \|\xtil(s)\|^2 + \tfrac{4 c'\beta_R}{\alpha_R} \epsilon^2.
\end{align}
We recognize $2c\gamma \le \calC_1$, and $\tfrac{4 c' \beta_R}{\alpha_R} = \calC_2 - \frac{1}{2\alpha_R}$, and invoke \Cref{lem:comparison} to obtain the desired bound.
\end{proof}
\begin{proof}[\textbf{Proof of \Cref{thm:RHC}}]
We first assume that $\alpha_V>0$ and bound the rate of covegence in terms of $\calC_0$; the steps of the proof can be repeated by replacing $\calC_0$ with $\calC_0^\delta$ throughout and then applying the tighter of the two bounds. 

For each $k \in \N$ we will let $\tilde{J}_{k} = J_T(\tilde{u}_{k}(\cdot;x_0,\uwarm_0),x_0)$ and $\tilde{J}_{k} = J_T(\tilde{u}_{k}(\cdot;x_0,\uwarm_0),x_0)$. Respectively, these are the cost incurred by the $k$-th planning solution and the $k$-th warm-start initial guess. For simplicity, we drop the dependence on $x_0$ and $\bar{u}_0$ from here on. We also let $(\tilde{x}_k(\cdot),\tilde{u}(\cdot))$ and $(\bar{x}_k(\cdot),\bar{u}_k(\cdot))$ denote the corresponding system trajectories. 
Our goal is to show that the sequence of losses $\set{\tilde{J}_k}_{k=1}^{\infty}$ is geometrically decreasing. Indeed, using property $1)$ of \Cref{asm:planner} we have

\begin{align}\label{eq:tot_bound}
\tilde{J}_{T}^{k+1} &\leq \bar{J}_{T}^{k+1} = \bar{J}_{T}^{k+1} - \tilde{J}_T^{k} + \tilde{J}_T^k\\\nonumber
&=\int_{T-\delta}^{T}Q(\bar{x}_{k+1}(\tau)) +R(\bar{u}_{k+1}(\tau))d\tau + V(\bar{x}_{k+1}(T))\\ \nonumber
&-\int_{0}^\delta Q(\tilde{x}_k(\tau)) +R(\tilde{u}_k(\tau))d \tau - V(\tilde{x}_k(T)) + \tilde{J}_{T}^k\\ \nonumber
&\leq \int_{T-\delta}^{T}\alpha_{Q}\|\bar{x}_{k+1}(\tau)\|_2^2  + \alpha_{V}\|\bar{x}_{k+1}(T))\|_2^2 \\ \nonumber
&+  \left( e^{-\frac{\delta}{\calC_1}} + \frac{\alpha_Q}{2}\calC_2  \epsilon_0^2\right)\tilde{J}_k
\end{align}
where the second inequality follows from \Cref{asm:cost_smooth}, the fact that $\bar{u}_{k+1}(t) = 0$ for each $t \in \sp{T-\delta, T}$, and applying the second inequality in \eqref{eq:exp_value_bound} which shows that $J_{T-\delta}(\tilde{u}_k|_{\sp{\delta,T}}, \tilde{x}_k(\delta)) \leq e^{-\frac{\delta}{\calC_1}}\tilde{J}_k + \frac{\alpha_Q}{2} \calC_2 \epsilon_0^2 \tilde{J}_k)$, where we have also used property $2)$ of \Cref{asm:planner}. The second fact also implies that $|\dot{\bar{x}}_{k+1}(t)| = |F(x_{k+1}(t), 0)| \leq L_F|\bar{x}(\bar{x}(t))|$ (by \Cref{ass:v_field_regularity}) for each $t\in \sp{T-\delta, T}$. Thus, by a standard application of a Gronwall-type inequality for each $t \in \sp{T-\delta, T}$ we have will will have 
\begin{align}\label{eq:term_bound}
\|\bar{x}_{k+1}(t)\|_2^2 &\leq e^{2L_F \delta} \|\bar{x}_{k+1}(T-\delta)\|_2^2\\
&\leq e^{2L_F\delta} \calC_0\left(\calC_1e^{-\frac{T-\delta}{\calC_1}}\tilde{J}_k + \frac{1}{2}\calC_2\epsilon_0^
2 \tilde{J}_k\right)
\end{align}
where we have used the fact that $\tilde{x}_{k}(T-\delta) = \bar{x}_{k+1}(T-\delta)$, property $2$ from \Cref{asm:planner} and \Cref{thm:plan}. Combining the above observations with the final inequality in \eqref{eq:tot_bound},integrating, and rearranging terms, and simplifying provides
\begin{multline}
\tilde{J}_{k+1} \leq \bigg(e^{-\frac{\delta}{\calC_1}} + \calC_0\bigg[\calC_1e^{2L_F\delta}(\delta \alpha_{Q} + \alpha_V)e^{-\frac{T-\delta}{\calC_1}}   + \frac{1}{2}e^{2L_F\delta}((\delta +1)\alpha_Q + \alpha_V) \calC_2  \epsilon_0^2 \bigg]\bigg)\tilde{J}_k
\end{multline}
Which simplifies to 
\begin{equation}\label{eq:geometric_decay}
\tilde{J}_{k+1} \leq \left(e^{-\frac{\delta}{\calC_1}} + \left[\mathcal{T}(\delta,T) + \mathcal{E}(\delta, \epsilon_0)\right]\right)\tilde{J}_{k}.
\end{equation}
 This geometric decay implies that 
\begin{equation}
J_T(\tilde{u}_k(0), \tilde{x}_k(0)) \leq  e^{ 2\eta(\delta, T,\epsilon_0) t_k} J_T(\tilde{u}_0(0), \tilde{x}_0(0))
\end{equation}
This can then be converted into the desired bound on the state trajectory by applying \Cref{lem:underbound,lem:overbound}.
\end{proof}
 }
{ 

In view of \Cref{lem:comparison}, to obtain a bound on $\Vtil(s)$ it suffices to bound  $\Vjacst(s) :=\inf_{\ubar_{[s,T]}} \Jlin_{T-s}(\cdot, ~ \xtil(s), \util_{[s,T]})$. Moreover, we can bound $\Vjacst(s)$ by bounding $\Jlin_{T-s}(\ubar_{[s,T]};~ \xtil(t), \util_{[s,T]})$ for any (possibly suboptimal) control $\ubar_{[s,T]}$; for simplicity, let us drop the ${[s,T]}$-subscript going forward. We select $\ubar(t) = \ubar_1(t) + \ubar_2(t)$, where $\ubar_1(t)$ satisfies $\Btil(t) \ubar_1(t) =- \tilde{d}(t)$, and where $\ubar_2$ witnesses $\gamma$-stabilizability at time $s$ as in \Cref{asm:stabilizable}.

With this choice of $\ubar(t)$ the dynamics of $\xbar(t)$ in $\Jlin_{T-s}$  are $\tfrac{\rmd}{\rmd t} \xbar(t) = \Atil(t)\xbar(t) + \Btil(t)\ubar(t) + \tilde{d}(t) = \Atil(t)\xbar(t) + \Btil \bar{u}_2(t)$, and, writing out $\Jlin_{T-s}$ explicitly, we obtain
\begin{align} 
\Vjacst(s) \le \textstyle \int_{t=s}^T (Q(\xbar(t)) + R(\ubar(t))\rmd t + V(\xbar(T)).
\end{align}
By the elementary bound $\|\ubar(t)\|^2 \le 2\|\ubar_1(t)\|^2 + 2\|\bar{u}_2(t)\|^2$, the following holds for constant $c = \max\{\beta_V,2\beta_R,\beta_Q\}$,
\begin{align}\textstyle
&\Vjacst(s)\textstyle\le   2\beta_R \int_{t=s}^T \|\ubar_1(t)\|^2 \dt \label{eq:P_line_one}\\
&\quad+c \left(\textstyle\int_{t=s}^T (\|\xbar(t)\|^2 + \|\ubar(t)\|^2)\rmd t + \|\xbar(T)\|^2\right). \label{eq:P_line_two}
\end{align}
To bound \eqref{eq:P_line_two}, we observe that $\ddt \xbar(t) = \Atil(t)\xbar(t) + \Btil(t)\ubar(t) + \tilde{d}(t) = \Atil(t)\xbar(t) + \Btil(t)\bar{u}_2(t)$, which corresponds to the $\hat{x}(t)$ dynamics in the definion of $\gamma$-stabilizability; thus, \eqref{eq:P_line_two} is at most $c \cdot \gamma \|\xtil(s)\|^2$. 

To bound \eqref{eq:P_line_one}, we use \eqref{asm:linearizable_bound} to bound $ \|\ubar_1(t)\|^2 \le  2L_x^2\|\xtil(t)\|^2 + 2 L_u^2\|\util(t)\|^2 \le 2c'(Q(\xtil(t)) +  R(\util(t))/\beta_R$, where $c' =  \beta_R\max\{\frac{L_x^2}{\alpha_Q},\frac{ L_u^2}{\alpha_R})$. Hence,  in view of \Cref{lem:comparison}, and the principle of optimality:
\begin{align*}
& 2 \beta_R\textstyle\int_{t=s}^T \|\bar{u}_1(t)\|^2 \dt \le  4c'\beta_R \textstyle\int_{t=s}^T (Q(\xtil(t)) +  R(\util(t))\dt \\
&\quad \le  4c' \beta_R \Vtil(s)  \le 4c'\beta_R(\Vjacst(s)  + \tfrac{\epsilon^2}{2\alpha_R})
\end{align*}
Putting the bounds together and rearranging:
\begin{align*}
(1-   4 \beta_R c' )\Vjacst(s) \le c \cdot \gamma \|\xtil(s)\|^2 + \tfrac{2 c'\beta_R}{ \alpha_R} \epsilon^2
\end{align*}
Under \Cref{asm:linearizable_bound}, we have $4\beta_R c'\le 1/2$, so that
\begin{align}
\Vjacst(s) \le 2c\cdot \gamma \|\xtil(s)\|^2 + 4 c'\beta_R \epsilon^2/\alpha_R.
\end{align}
We recognize $2c\gamma \le \calC_1$, and $\tfrac{4 c' \beta_R}{\alpha_R} = \calC_2 - \frac{1}{2\alpha_R}$, and invoke \Cref{lem:comparison} to obtain the desired bound. \qedd

}
\bibliographystyle{IEEEtran}
\bibliography{IEEEabrv,refs.bib}

\end{document}